# Modular forms on the moduli space of polarised K3 surfaces

Arie Peterson[*]

**Abstract**

We study the moduli space $\mathcal{F}_{2d}$ of polarised K3 surfaces of degree $2d$. We compute all relations between Noether–Lefschetz divisors on these moduli spaces for $d$ up to around 50. This leads to a very concrete description of the rational Picard group of $\mathcal{F}_{2d}$. We show how to determine the coefficients of boundary components of relations in the rational Picard group, giving relations on a (toroidal) compactification of $\mathcal{F}_{2d}$. We draw conclusions from this about the Kodaira dimension of $\mathcal{F}_{2d}$, in many cases confirming earlier results by Gritsenko, Hulek and Sankaran [12], and in two cases giving a new result.

This is an abridged version of the PhD thesis [28].

## 1 Introduction

We consider the moduli space $\mathcal{F}_{2d}$ of polarised K3 surfaces of degree $2d$. For every positive integer $d$ the moduli space $\mathcal{F}_{2d}$ is a 19-dimensional irreducible quasi-projective variety, by the work of Baily–Borel [2]. There is a period map giving an alternative description of $\mathcal{F}_{2d}$ as an arithmetic quotient (a quotient of a complex domain by a discrete group) associated to the lattice

$$L_{2d} = \langle -2d \rangle \oplus U^{\oplus 2} \oplus E_8(-1)^{\oplus 2} , \tag{1}$$

where $U$ is the hyperbolic plane and $E_8$ is the root lattice of the Lie algebra of the same name. In terms of this lattice the 19-dimensional complex domain $\mathcal{D}_{2d}$ is given by

$$\mathcal{D}_{2d} \cup \overline{\mathcal{D}_{2d}} = \{\mathbb{C}z : (z,z) = 0, (z,\bar{z}) > 0\} \subset \mathbb{P}(L_{2d} \otimes \mathbb{C}) \tag{2}$$

and $\mathcal{F}_{2d}$ is isomorphic to the arithmetic quotient $\tilde{\mathrm{O}}^+(L_{2d}) \backslash \mathcal{D}_{2d}$, where $\tilde{\mathrm{O}}^+(L_{2d})$ is the group of stable automorphisms of spinor norm 1 of the lattice $L_{2d}$. The fact that the period map indeed gives an isomorphism from $\mathcal{F}_{2d}$ to this quotient was proved by Piatetskii-Shapiro and Shafarevich [30] and Todorov [33].

An important feature of $\mathcal{F}_{2d}$ is its Picard group. A result by O'Grady [27] already showed that as the polarisation degree $2d$ increases, the rank of $\mathrm{Pic}(\mathcal{F}_{2d})$ is unbounded. The theory of modular forms can be used to study this Picard group; this was accelerated by the work of Borcherds [4], allowing the construction of many modular forms on arithmetic quotients. In particular, Bruinier [6]

---

[*]This research was funded by NWO grant 613.001.018.

computed the dimension of the part of the rational Picard group generated by
Noether–Lefschetz divisors, and recently it was shown [3] that this part in fact
equals the full rational Picard group.

In sections 3.3 and 3.4 we show how to compute a complete and explicit
description of the rational Picard group of $\mathcal{F}_{2d}$. We have used this procedure
to get an explicit basis of $\mathrm{Pic}_{\mathbb{Q}}(\mathcal{F}_{2d})$ for $d$ up to 50 and coefficients expressing
many Noether–Lefschetz divisors $H(\gamma, n)$ in terms in this basis.

Also, we show how to complete any relation in $\mathrm{Pic}_{\mathbb{Q}}(\mathcal{F}_{2d})$ to the boundary
of a specific toroidal compactification $\overline{\mathcal{F}_{2d}}$, getting relations in $\mathrm{Pic}_{\mathbb{Q}}(\overline{\mathcal{F}_{2d}})$; see
theorem 4.8. Additionally, we introduce a new method of "theta ghosts" to
compute bounds for the variation of the coefficients of a relation in $\mathrm{Pic}_{\mathbb{Q}}(\overline{\mathcal{F}_{2d}})$
over the (very large) set of 1-cusps; see section 4.5.

Having a good description of the Picard group $\mathrm{Pic}_{\mathbb{Q}}(\mathcal{F}_{2d})$, we would like to
determine the structure of the effective cone inside it. It is not obvious though
whether the effective cone is generated by classes of Noether–Lefschetz divisors;
this seems to be an interesting geometric question. We have looked in some
detail at the subcone generated by Noether–Lefschetz classes: see section 3.5.

Another important question about $\mathcal{F}_{2d}$ concerns its Kodaira dimension. For
$d \in \{1, 2, 3, 4\}$ there are well-known explicit constructions of polarised K3 surfaces and $\mathcal{F}_{2d}$ is unirational in those cases. Mukai has extended this to some
slightly higher values of $d$, proving that $\mathcal{F}_{2d}$ is unirational (hence $\kappa(\mathcal{F}_{2d}) = -\infty$)
also for $d \in \{5, 6, 7, 8, 9, 10, 11, 12, 15, 17, 19\}$ (see [20], [22], [23], [24], [21]). In
the other direction, Kondō has shown in [17] that $\mathcal{F}_{2p^2}$ is of general type for
large enough primes $p$, and later in [18] that the Kodaira dimension of $\mathcal{F}_{2d}$ is
non-negative if

$$d \in \{42, 43, 51, 53, 55, 57, 59, 61, 66, 67, 69, 74, 83, 85, 105, 119, 133\} \ .$$

Gritsenko, Hulek and Sankaran more recently proved [12] that $\mathcal{F}_{2d}$ is of general
type for all $d > 61$ and also for $d \in \{46, 50, 54, 57, 58, 60\}$, and that its Kodaira
dimension is non-negative for $d \in \{40, 42, 43, 46, 48, 49, 51, 52, 53, 55, 56, 59, 61\}$.
It was noted by the author and Sankaran in [29] that their method also applies
to the case $d = 52$, proving that $\mathcal{F}_{2 \cdot 52}$ is of general type.

We have confirmed and extended these results. Firstly, we have some concrete additions for low $d$:

**Theorem** (3.21). *If $1 \leq d \leq 15$, the moduli space $\mathcal{F}_{2d}$ has Kodaira dimension $-\infty$.*

This gives the new cases $d \in \{13, 14\}$ and reproves Mukai's results for $1 \leq d \leq 12$
and $d = 15$ in a completely different way.

Secondly, we have results for cases of intermediate $d$, some of which are
conditional:

**Theorem** (3.24). *Let $16 \leq d \leq 39$ or $d \in \{41, 44, 45, 47\}$. Either $\kappa(\mathcal{F}_{2d}) = -\infty$
or there exists an irreducible codimension 1 subvariety of $\mathcal{F}_{2d}$ that is not a
Noether–Lefschetz divisor.*

**Theorem** (4.14,3.26). *If $d \in \{40, 42, 43, 48, 49, 55, 56\}$, then $\kappa(\mathcal{F}_{2d}) \geq 0$. Also,
either $\kappa(\mathcal{F}_{2d}) < 19$ or there exists an irreducible codimension 1 subvariety of $\mathcal{F}_{2d}$
that is not a Noether–Lefschetz divisor.*

Thirdly, we have reproved some of the results of [12], where $d$ is slightly higher: if $d \in \{46, 50, 52, 54\}$, then $\mathcal{F}_{2d}$ is of general type; see theorem 4.12. Also, for $d \in \{40, 42, 43, 48, 49, 55, 56\}$ the Kodaira dimension of $\mathcal{F}_{2d}$ is non-negative; see theorem 4.14.

Our computations are very explicit: in the cases where we find that $\kappa(\mathcal{F}_{2d})$ is non-negative, we get a concrete linear relation for the canonical divisor on the open part of the moduli space in terms of Noether–Lefschetz divisors.

## 2 Background and definitions

A K3 surface is a smooth surface $S$ such that $\omega_S \cong \mathcal{O}_S$ and $\mathrm{H}^1(S, \mathcal{O}_S) = 0$. We will only consider projective K3 surfaces over $\mathbb{C}$ in this paper. If $S$ is a K3 surface, then the rank of its second cohomology group $\mathrm{H}^2(S, \mathbb{Z})$ is 22. Moreover, the lattice structure on $\mathrm{H}^2(S, \mathbb{Z})$ (given by the cup product) is completely determined: $\mathrm{H}^2(S, \mathbb{Z})$ is isomorphic to the K3 lattice $L_{K3} = U^{\oplus 3} \oplus E_8(-1)^{\oplus 2}$.

We will be interested in the moduli space of polarised K3 surfaces. For us, a polarisation is a choice of primitive nef line bundle $H$ with positive self-intersection $H^2 = H \cdot H$. (This is called a quasi-polarisation by some, who then require a polarisation to be ample.) The so-called degree $H^2$ is always even, and we will write it as $H^2 = 2d$ for a positive integer $d$. (Note that this implies that $g = d + 1$ is the genus of a smooth section of the polarisation class $H$.)

A lattice is a free $\mathbb{Z}$-module $L$ of finite rank together with a symmetric bilinear map $(\cdot, \cdot) : L \times L \to \mathbb{Q}$. As general references for the theory of lattices, see [10],[25],[11],[16].

If $L$ is an integral lattice, we denote the discriminant group $L^\vee / L$ by $D_L$, and the group of automorphisms of $L$ by $\mathrm{O}(L)$. The subgroup $\tilde{\mathrm{O}}(L) \subseteq \mathrm{O}(L)$ is the kernel of the map $\mathrm{O}(L) \to \mathrm{O}(D_L)$; in other words, it is the set of automorphisms of $L$ that act trivially on the discriminant group. The group $\mathrm{O}^+(L) \subseteq \mathrm{O}(L)$ is the subgroup of automorphisms of spinor norm 1 (the spinor norm is multiplicative, and the spinor norm of a reflection in a vector of positive norm is $-1$). The group $\tilde{\mathrm{O}}^+(L) \subseteq \mathrm{O}(L)$ is the intersection of $\tilde{\mathrm{O}}(L)$ and $\mathrm{O}^+(L)$.

If $L$ is a lattice, we define the *genus* of $L$ (denoted by $\mathcal{G}(L)$) to be the set of isomorphism classes of lattices $K$ that are locally isomorphic to $L$, i.e. such that for all primes $p$ (including the infinite prime) we have $K \otimes \mathbb{Z}_p \cong L \otimes \mathbb{Z}_p$ (with the convention that $\mathbb{Z}_\infty = \mathbb{R}$).

If we restrict to even lattices, the local equivalence class of a lattice at all finite primes is captured precisely by the discriminant module:

**Theorem 2.1** ([25, 1.9.4]). *Even lattices are locally equivalent if and only if they have the same signature and isomorphic discriminant modules.*

### 2.1 Vector-valued modular forms associated to a lattice

Throughout this section, let $L$ be an even lattice and $k \in \frac{1}{2}\mathbb{Z}$. Following the conventions of [6], we introduce a class of modular forms associated to $L$ of weight $k$ with respect to the metaplectic group, taking values in the vector space $\mathbb{C}[D_L]$. First of all, we should specify how the metaplectic group acts on vector-valued functions on the upper half-plane.

**Definition 2.2.** Let $f : \mathbb{H} \to \mathbb{C}[D_L]$ and $(A, \varepsilon) \in \mathrm{Mp}_2(\mathbb{Z})$. The slash operator $|_k^*$ is defined by
$$f |_k^* (A, \varepsilon)(z) = \varepsilon(z)^{-2k} \rho_L^\vee (A, \varepsilon)^{-1} f(Az) \ . \tag{3}$$
Here, $\rho_L^\vee$ is the dual of the Weil representation $\rho_L$ of $\mathrm{Mp}_2(\mathbb{Z})$ on $\mathbb{C}[D_L]$, defined in section 3.2 of [28].

**Definition 2.3.** A holomorphic function $f : \mathbb{H} \to \mathbb{C}[D_L]$ is called a modular form of weight $k$ on the metaplectic group $\mathrm{Mp}_2(\mathbb{Z})$ if $f$ is holomorphic at infinity and for every $(A, \varepsilon) \in \mathrm{Mp}_2(\mathbb{Z})$ we have $f |_k^* (A, \varepsilon) = f$.

Suppose that $\varphi \in M(k, L)$ is a vector-valued modular form. We may write such a form as a Fourier expansion
$$\varphi(\tau) = \sum_{\gamma, n} a_{\gamma, n} q^n \mathbf{e}_\gamma \ , \tag{4}$$
where we write $q = \mathrm{e}(\tau) = e^{2\pi i \tau}$, and $\mathbf{e}_\gamma$ is the standard basis vector of $\mathbb{C}[D_L]$ associated to $\gamma \in D_L$.

**Definition 2.4.** Given $(\gamma, n)$ (such that $n \in -\gamma^2/2 + \mathbb{Z}$ and $n \geq 0$), let $c_{\gamma,n} : M(k, L) \to \mathbb{C}$ be the function taking a form $\varphi$ to its $(\gamma, n)$-coefficient $a_{\gamma, n}$.

The modular curve $\mathrm{SL}_2(\mathbb{Z}) \backslash \mathbb{H}$ on which these vector-valued modular forms live has one cusp. However, we may want to distinguish the behaviour of the different components of a vector-valued form at the cusp. This leads to the following useful abuse of language:

**Definition 2.5.** Let $L$ be an even lattice. The cusps of $L$ are the isotropic elements of the discriminant group of $L$ (i.e., the elements $\gamma \in D_L = L^\vee/L$ such that $\gamma^2/2 = 0 \in \mathbb{Q}/\mathbb{Z}$).

Using this terminology, we might say that a given vector-valued modular form $\varphi$ vanishes at the cusp $\gamma \in D_L$ if the $\mathbf{e}_\gamma$-component of $\varphi$ vanishes at the single cusp of the modular curve.

**Remark 2.6.** There is a tight connection between the vector-valued modular forms associated to $L$ and modular forms on the arithmetic quotient associated to $L$. The cusps of $L$ in the above sense correspond with the 0-dimensional cusps of that arithmetic quotient, giving an alternative interpretation of definition 2.5.

**Definition 2.7.** We say that a modular form $\varphi \in M(k, L)$ is a cusp form if its coefficients satisfy $c_{\gamma, 0}(\varphi) = 0$ for all cusps $\gamma$ of $L$.

As a non-standard extension of this terminology, we say that $\varphi$ is an almost cusp form if $c_{\gamma, 0}(\varphi) = 0$ for all isotropic $\gamma \in D_L$ except perhaps the standard cusp $\gamma = \bar{0} \in D_L$.

We denote the space of cusp forms of weight $k$ by $S(k, L)$ and the space of almost cusp forms of weight $k$ by $AC(k, L)$.

Note that we have $S(k, L) \subseteq AC(k, L) \subseteq M(k, L)$. If $k > 2$, the first inclusion is strict, and the last inclusion is an equality if and only if $L$ has only one cusp: this follows from the existence of an Eisenstein series $E_\gamma$ at every cusp $\gamma$ of $L$ (see [9] and section 3.4 of [28]).

### 2.1.1 Vector-valued theta series

Let $L$ be a positive-definite even lattice. The theta series $\Theta_L$ counts the number of vectors in $L$ of all possible lengths:

$$\Theta_L = \sum_{v \in L} q^{v^2/2} \ . \tag{5}$$

If $L$ is non-unimodular, we may collect more information by considering vectors in the dual lattice $L^\vee$. To keep track of the discriminant class (i.e., the coset of the vector in $D_L = L^\vee/L$) we let the theta series take values in the group algebra $\mathbb{C}[D_L]$:

$$\boldsymbol{\Theta}_L = \sum_{v \in L^\vee} q^{v^2/2} \mathbf{e}_{v+L} \ . \tag{6}$$

We call this the vector-valued theta series associated to the lattice $L$.

**Proposition 2.8** ([4, section 4]). *If $L$ is a definite even lattice of rank $k$, then $\boldsymbol{\Theta}_L$ is a vector-valued modular form (see section 2.1) of weight $k/2$, with values in the Weil representation $\rho_L$ (see section 3.2 of [28]).*

## 3 The moduli space of polarised K3 surfaces: open part

### 3.1 Heegner and Noether–Lefschetz divisors

The two descriptions of $\mathcal{F}_{2d}$ (as a locally symmetric variety and as a moduli space) give rise to two sets of natural divisors which turn out to coincide. We will first introduce the divisors on the arithmetic side.

**Definition 3.1.** Given a vector $v \in L^\vee$, the Heegner divisor $H_v$ is the subset of $\mathcal{D}_{2d} \subset \mathbb{P}(L \otimes \mathbb{C})$ orthogonal to $v$.

Note that this is a divisor on the symmetric domain $\mathcal{D}_{2d}$, not on the arithmetic quotient $\mathcal{F}_{2d}$. As any multiple of $v$ will have the same orthogonal complement, we may as well assume $v$ to be primitive in $L^\vee$.

Now, the orbits of primitive vectors of $L^\vee$ under the action of $\tilde{O}^+(L)$ are classified exactly by the coset $v + L \in D_L$ and the square $v^2$ ([15, lemma 7.5]). This motivates the following definition.

**Definition 3.2.** Let $\gamma \in D_L$ and $n \in \mathbb{Z} - \gamma^2/2$ such that $n < 0$. We denote the corresponding Heegner divisor by $H(\gamma, n)$:

$$H(\gamma, n) = \tilde{O}^+(L) \Big\backslash \Big( \sum_{\substack{v \in L^\vee \\ v \equiv \gamma \bmod L \\ v^2 = 2n}} H_v \Big) \ . \tag{7}$$

Note that $H(\gamma, n)$ has multiplicity 1 everywhere if $-\gamma \neq \gamma$, but it has multiplicity 2 everywhere if $-\gamma = \gamma$.

Apart from this, the tautological bundle $\mathcal{O}(-1)$ on $\mathbb{P}(L \otimes \mathbb{C})$ descends to a line bundle $\lambda$ on $\mathcal{F}_{2d}$, the class of which we will somewhat loosely refer to as the Hodge class. (Note that $\lambda$ is isomorphic to $\pi_*\omega_\pi$, the pushforward of the relative dualising sheaf of the universal K3 surface $\pi: \mathcal{X}_{2d} \to \mathcal{F}_{2d}$.) The divisor (pick any) associated to the dual bundle $-\lambda$ is often considered together with the Heegner divisors, and denoted by $H(\overline{0}, 0)$. Note the minus sign: while the Heegner divisors are effective, this $H(\overline{0}, 0)$ is anti-ample. This may seem a strange choice, but it turns out to be natural.

We now introduce the so-called Noether–Lefschetz divisors on $\mathcal{F}_{2d}$: these are the loci where the Picard group of the K3 surface jumps from rank 1 to rank 2. Because the rank 2 Picard group can be one of infinitely many non-isomorphic lattices, we get an infinite set of divisors. There are two slightly different ways to proceed.

One option is to prescribe the isomorphism class of the Picard lattice (the Picard group with its intersection form). This will give irreducible divisors; we will look at these in section 3.2. Right now, we take a slightly different approach: we look at the locus of polarised K3 surfaces that have an extra class in their Picard group with given intersection numbers.

**Definition 3.3.** Given $h \in \mathbb{N}, a \in \mathbb{Z}$ such that $a^2 - 4d(h-1)$ is positive, the Noether–Lefschetz divisor $D_{h,a} \subset \mathcal{F}_{2d}$ is supported on the locus of polarised K3 surfaces $(S, H)$ that have a divisor class $\beta \in \operatorname{Pic} S$ (with $\beta$ not in the span of $H$) of square $\beta^2 = 2h - 2$ and degree $\beta \cdot H = a$.

The multiplicity of the irreducible component of $D_{h,a}$ consisting (generically) of K3 surfaces that have Picard lattice isomorphic to the rank 2 lattice $L$ by definition equals the number of elements $\beta \in L$ with $\beta^2 = 2h-2$ and $\beta \cdot H = a$ (where we write $H$ for the element of $L$ corresponding to the polarisation class).

**Remark 3.4.** The irreducible components of $D_{h,a}$ are exactly the so-called irreducible Noether–Lefschetz divisors, each of which parametrises K3 surfaces with rank 2 Picard lattice of a fixed isomorphism class. We will look at those in section 3.2.

As we alluded to before, these geometrically defined Noether–Lefschetz divisors coincide with the Heegner divisors:

**Lemma 3.5** ([19, section 4.4, lemma 3])**.** *Suppose that $\gamma = aw/2d \in D_L$ and $n = h - 1 - a^2/4d$. Then we have $D_{h,a} = H(\gamma, n)$.*

Therefore, we will use the terms "Heegner divisor" and "Noether–Lefschetz divisor" interchangeably.

Having introduced a large set of divisors on the moduli space $\mathcal{F}_{2d}$, we would like to understand their role in the Picard group $\operatorname{Pic}_{\mathbb{Q}}(\mathcal{F}_{2d})$.

**Definition 3.6.** The Noether–Lefschetz Picard group $\operatorname{Pic}_{\mathbb{Q}}^{\mathrm{NL}}(\mathcal{F}_{2d})$ is the $\mathbb{Q}$-subspace of $\operatorname{Pic}_{\mathbb{Q}}(\mathcal{F}_{2d})$ spanned by $H(\overline{0}, 0)$ and the (infinite) set of Heegner divisors $H(\gamma, n)$, for all $\gamma \in D_L$ and all $n \in \mathbb{Z} - \gamma^2/2$ with $n < 0$.

A first natural question is, whether the subspace $\operatorname{Pic}_{\mathbb{Q}}^{\mathrm{NL}}(\mathcal{F}_{2d})$ might in fact be equal to the full Picard group $\operatorname{Pic}_{\mathbb{Q}}(\mathcal{F}_{2d})$. Recently, this has been proved by Bergeron, Li, Millson and Moeglin:

**Theorem 3.7** ([3])**.** *The rational Picard group is spanned by Noether–Lefschetz divisors, so indeed $\operatorname{Pic}_{\mathbb{Q}}^{\mathrm{NL}}(\mathcal{F}_{2d}) = \operatorname{Pic}_{\mathbb{Q}}(\mathcal{F}_{2d})$.*

The proof uses the representation theory of the orthogonal groups occurring in the arithmetic description of $\mathcal{F}_{2d}$ as a locally symmetric domain.

A next question is to compute the dimension of $\mathrm{Pic}_{\mathbb{Q}}(\mathcal{F}_{2d})$ as a function of $d$. A formula for $\dim \mathrm{Pic}_{\mathbb{Q}}^{\mathrm{NL}}(\mathcal{F}_{2d})$ has been found by Bruinier; see [7] and [28, Theorem 4.2.10]. By the above theorem by Bergeron, Li, Millson and Moeglin, Bruinier's formula in fact gives the dimension of $\mathrm{Pic}_{\mathbb{Q}}(\mathcal{F}_{2d})$.

## 3.2 Irreducible Noether–Lefschetz divisors

Let $(S, H)$ be a polarised K3 surface. The moduli point in $\mathcal{F}_{2d}$ corresponding to this surface lies on the Noether–Lefschetz divisor $D_{h,a}$ if and only if there exists a divisor class on $S$ with intersection numbers $h, a$. This forces the Picard lattice of $S$ to be at least rank 2. Even if we suppose that the Picard lattice of $S$ has rank exactly 2, we cannot determine the lattice structure of $\mathrm{Pic}(S)$ from this condition alone. The reason is that in general there are several non-isomorphic rank 2 lattices having an element with the prescribed intersection numbers.

This is relevant from a geometric point of view: because of this phenomenon, the Noether–Lefschetz divisors $D_{h,a}$ are not irreducible. The irreducible components are the loci where the Picard lattice of the polarised K3 surface is of a given isomorphism class. We first present a convenient parametrisation of these isomorphism classes.

**Definition 3.8.** Let $(L, H)$ be a $2d$-polarised even lattice of rank 2 and signature $(1, 1)$. Then the discriminant $\Delta \in \mathbb{Z}_+$ and coset $\delta \in (\mathbb{Z}/2d\mathbb{Z})/\pm$ of $L$ are defined as follows: first choose an element $\Gamma \in L$ such that $H$ and $\Gamma$ form a basis of $L$. We write down the intersection matrix of $L$ with respect to this basis:

$$\begin{pmatrix} 2d & y \\ y & 2x \end{pmatrix}$$

(so $y = H \cdot \Gamma$ and $2x = \Gamma^2$). Then $\Delta = y^2 - 4dx$, and $\delta = \overline{y} \in (\mathbb{Z}/2d\mathbb{Z})/\pm$.

(Note that the number $\Delta$ is in fact minus the usual discriminant of a lattice (i.e., the determinant of the intersection matrix); we prefer to work with positive numbers.) It is an easy exercise to check that the discriminant and coset are independent of the choice of element $\Gamma$, so they are invariants of this particular class of polarised lattices. In fact, they are complete invariants:

**Proposition 3.9** ([28, lemma 4.2.12, proposition 4.2.13]). *The $2d$-polarised even lattices of rank 2 and signature $(1, 1)$ are classified by their discriminant and coset. There exists a rank 2 even hyperbolic $2d$-polarised lattice of discriminant $\Delta$ and coset $\delta$ if and only if*

$$\Delta \equiv \delta^2 \mod 4d \ . \tag{8}$$

Now, using this description of rank 2 lattices in terms of their invariants $\Delta$ and $\delta$ we arrive at the following definition of the irreducible Noether–Lefschetz divisors.

**Definition 3.10.** Given $\Delta \in \mathbb{Z}_+$ and $\delta \in (\mathbb{Z}/2d\mathbb{Z})/\pm$, the Noether–Lefschetz divisor $P_{\Delta, \delta}$ is the closure of the set of polarised K3 surfaces $(S, H)$ such that $\mathrm{Pic}\, S$ has discriminant $\Delta$ and coset $\delta$.

By taking the closure, we include surfaces with a Picard lattice that contains the given lattice, but has rank higher than 2. The thesis [26] contains a proof that these divisors $P_{\Delta,\delta}$ are indeed irreducible if they are non-empty (which happens if and only if there exists a lattice with discriminant $\Delta$ and coset $\delta$).

Section 4.2.2 of [28] contains a formula that expresses which irreducibles $P_{\Delta,\delta}$ occur in a given $D_{h,a}$ and with what multiplicity, giving an easy way to express reducible Noether–Lefschetz divisors in terms of irreducible ones and vice versa.

### 3.3 Borcherds' construction of modular forms on $\mathcal{F}_{2d}$

In [4], Borcherds gives a construction of modular forms on the arithmetic quotient associated to a lattice of signature $(2,n)$. We apply this to our situation, where the lattice is $L = L_{2d}$ (of signature $(2,19)$), and the associated arithmetic quotient is exactly our moduli space $\mathcal{F}_{2d}$. This gives us a large supply of modular forms on our space $\mathcal{F}_{2d}$.

We introduce a few tools to describe the pole behaviour of vector-valued modular forms.

**Definition 3.11.** We let $\mathrm{Sing}(L)$ be the space of Laurent polynomials (in $q$) with values in $\mathbb{C}[D_L]$ having only non-positive powers of $q$.

Similarly, $\mathrm{Sing}(L)^-$ is the subspace of $\mathrm{Sing}(L)$ of elements having only negative powers of $q$, and $\mathrm{Sing}(L)_{\bar{0}}^-$ is the subspace of $\mathrm{Sing}(L)$ of elements having only negative powers of $q$ except possibly a term $q^0 \mathbf{e}_{\bar{0}}$.

We write $\mathrm{Obstruct}(k, L)$ for the obstruction space to the existence of vector-valued modular forms of weight $k$ with given principal part. More formally, $\mathrm{Obstruct}(k, L)$ is the quotient of $\mathrm{Sing}(L)$ by the image of the map taking a meromorphic vector-valued modular form of weight $k$ to its principal part.

By an application of Serre duality, the space $M(k, L)$ is dual to the obstruction space $\mathrm{Obstruct}(2-k, L)$ of obstructions to a given element of $\mathrm{Sing}(L)$ being the principal part of a meromorphic vector-valued modular form of weight $2-k$ and representation $\rho_L$. The duality is realised by the residue map

$$M(k, L) \times \mathrm{Sing}(L) \to \mathbb{C}$$
$$(\varphi, f) \mapsto \mathrm{Res}(\frac{\varphi f}{q^{1/N}} \, \mathrm{d}q^{1/N}) \ . \tag{9}$$

**Proposition 3.12** ([5]). *The above duality identifies the coefficient function $c_{\gamma,n} \in M(k, L)^\vee$ with $[q^{-n}\mathbf{e}_\gamma] \in \mathrm{Obstruct}(2-k, L)$. A linear combination*

$$\sum_{\gamma, n} a_{\gamma,n} c_{\gamma,n} \in M(k, L)^\vee \tag{10}$$

*is zero on $M(k, L)$ if and only if the corresponding obstruction*

$$[\sum_{\gamma, n>0} a_{\gamma,n} q^{-n} \mathbf{e}_\gamma] \in \mathrm{Obstruct}(2-k, L) \tag{11}$$

*vanishes, i.e., if and only if this equals the principal part of some meromorphic vector-valued form of weight $2-k$ and representation $\rho_L$.*

**Remark 3.13.** The analogous statement holds for the space $S(k, L)$ of cusp forms if we restrict the principal parts to $\mathrm{Sing}^-(L)$, and for the space $AC(k, L)$ of almost cusp forms if we restrict the principal parts to $\mathrm{Sing}_{\bar{0}}^-(L)$.

We may now formulate the result of Borcherds' construction.

**Theorem 3.14** ([4, Theorem 13.3]). *Let $f = \sum_{\gamma, n \geq 0} a_{\gamma, n} q^{-n} \mathbf{e}_\gamma \in \mathrm{Sing}(L)$ be the principal part of a meromorphic modular form of weight $1 - b^-/2 = 1 - 19/2 = -17/2$ (i.e., $[f] = 0 \in \mathrm{Obstruct}(-17/2, L)$). Assume that all coefficients $a_{\gamma, n}$ are integral. Then there is a meromorphic modular form $\Psi$ (scalar valued) on $\mathcal{F}_{2d}$, of weight $a_{\bar{0}, 0}/2$, with divisor $1/2 \cdot \sum_{\gamma, n > 0} a_{\gamma, n} H(\gamma, n)$. Moreover, $\Psi$ has the following product expansion around the cusp associated to $z$, in the tube domain parametrisation:*

$$\Psi_z(Z_M) = C\,\mathrm{e}((Z_M, \rho_M)) \prod_{\substack{\lambda \in M^\vee \\ (\lambda, W_M) > 0}} \prod_{\substack{\delta \in D_L \\ \delta|_M = \lambda}} \left(1 - \mathrm{e}((\lambda, Z_M) + (\delta, z'))\right)^{a_{\lambda, \lambda^2/2}} . \tag{12}$$

*Here, the number $C$ is some nonzero constant; $W_M$ is a Weyl chamber (with respect to $f$) that has $z$ in its closure (this is the subset of the period domain on which the expansion will be valid); $\rho_M = \rho(M, W_M, f) \in M \otimes \mathbb{Q}$ is the corresponding Weyl vector; the notation $(\lambda, W_M) > 0$ means that $(\lambda, w) > 0$ for all $w \in W_M$.*

About the condition $(\lambda, W_M) > 0$: if $\lambda$ is such that $a_{\lambda, \lambda^2/2} \neq 0$ (and note that other $\lambda$ do not contribute!), then it suffices to check this for any single $w_0 \in W_M$.

Note that even though $f \in \mathrm{Sing}(L)$, only the terms of $f$ with $n < 0$ or $(\gamma, n) = (\bar{0}, 0)$ are used, so we might as well take $f \in \mathrm{Sing}_0^-(L)$. The theorem shows that if the combination

$$a_{\bar{0}, 0} \mathbf{e}_{\bar{0}} + \sum_{\gamma, n > 0} a_{\gamma, n} q^{-n} \mathbf{e}_\gamma \in \mathrm{Sing}_0^-(L) \tag{13}$$

vanishes in the obstruction space $\mathrm{Obstruct}(-17/2, L)$, then

$$a_{\bar{0}, 0} H(\bar{0}, 0) + \sum_{\gamma, n > 0} a_{\gamma, n} H(\gamma, n) \tag{14}$$

is linearly equivalent to the zero divisor on $\mathcal{F}_{2d}$. (Recall that a modular form of weight $k$ is a section of the line bundle $\lambda^{\otimes k} \sim -k H(\bar{0}, 0)$.) On the other hand, by proposition 3.12 and remark 3.13 we know that

$$\sum_{\gamma, n \geq 0} a_{\gamma, n} q^{-n} \mathbf{e}_\gamma \in \mathrm{Sing}_0^-(L) \tag{15}$$

vanishes in $\mathrm{Obstruct}(-17/2, L)$ if and only if the corresponding functional of coefficients of almost cusp forms of weight $2 - (-17/2) = 21/2$ vanishes, i.e., if

$$\sum_{\gamma, n \geq 0} a_{\gamma, n} c_{\gamma, n} = 0 \in AC(d)^\vee . \tag{16}$$

Therefore, we may produce relations in $\mathrm{Pic}_\mathbb{Q}(\mathcal{F}_{2d})$ among Heegner divisors and the Hodge class $\lambda$ by computing equalities among coefficients of vector-valued almost cusp forms.

Moreover, Bruinier [8, Theorem 1.2] shows that any meromorphic modular form on $\mathcal{F}_{2d}$ with divisor supported on Heegner divisors occurs as a result of Borcherds' construction. Translating to geometric terms, this means that all relations among Noether–Lefschetz divisors come from linear combinations of coefficients that vanish on all almost cusp forms.

Finally, by the recent work of [3], the rational Picard group of $\mathcal{F}_{2d}$ is generated by Noether–Lefschetz divisors, so we may summarise all the above as follows.

**Theorem 3.15.** *The rational Picard group of $\mathcal{F}_{2d}$ is isomorphic to the dual of the space of rational vector-valued almost cusp forms of weight* $21/2$. *This isomorphism* $\mathrm{Pic}_{\mathbb{Q}}(\mathcal{F}_{2d}) \to AC(d)_{\mathbb{Q}}^{\vee}$ *sends* $[H(\gamma, n)]$ *to the coefficient function* $c_{\gamma, n} : AC(d)_{\mathbb{Q}}^{\vee} \to \mathbb{Q}$; *as a special case,* $\lambda = -[H(\bar{0}, 0)]$ *is sent to* $-c_{\bar{0}, 0}$.

*Proof.* As described above, this is a direct combination of Borcherds' construction of forms on arithmetic quotients [4], Bruinier's converse theorem [8], and the result [3] by Bergeron, Li, Millson and Moeglin that $\mathrm{Pic}_{\mathbb{Q}}(\mathcal{F}_{2d})$ is generated by Noether–Lefschetz divisors. □

## 3.4 Computing relations in $\mathrm{Pic}_{\mathbb{Q}}(\mathcal{F}_{2d})$

As we saw in the previous section, the relations among divisors on $\mathcal{F}_{2d}$ are exactly given by linear relations between coefficients of vector-valued modular forms.

We have implemented the method suggested by Raum [31] to compute a basis of that space of vector-valued modular forms up to any wanted number of Fourier coefficients (see [28, section 3.6] more for details). These data then give all relations among Noether–Lefschetz divisors in concrete form, and by theorem 3.7, that gives us a complete description of the rational Picard group.

**Example 3.16.** Let us give an example of the data that result from this procedure in the case $d = 1$. Recall that the space $M(1)$ of vector-valued modular forms has dimension 2, and that there is only a single cusp, so that the space of almost cusp forms is the whole space $AC(1) = M(1)$.

As a basis for $AC(1)^{\vee}$, we get $\{\varphi_1, \varphi_2\}$, where

$$\varphi_1 = c_{\bar{0}, 0}, \quad \varphi_2 = \frac{105457575250}{169227} c_{\bar{0}, 0} + c_{\bar{0}, -1}. \tag{17}$$

Employing the isomorphism between $AC(d)^{\vee}$ and $\mathrm{Pic}_{\mathbb{Q}}(\mathcal{F}_{2d})$, we conclude that a basis of $\mathrm{Pic}_{\mathbb{Q}}(\mathcal{F}_{2 \cdot 1})$ is formed by

$$H(\bar{0}, 0) = -\lambda \quad \text{and} \quad -\frac{105457575250}{169227} \lambda + H(\bar{0}, -1). \tag{18}$$

Alternatively, we might take as a basis the class $\lambda$ and the class of the Noether–Lefschetz divisor $H(\bar{0}, -1) = D_{0, 0}$.

Given any other Noether–Lefschetz divisor $D_{h, a} = H(\gamma, n)$, we may read off its coefficients with respect to the basis $\{\varphi_1, \varphi_2\}$ directly from the corresponding coefficients of $\varphi_1$ and $\varphi_2$. For instance, take $D_{0,1} = H(\overline{1}, -1/4)$:

$$\begin{aligned}
D_{0,1} &= H(\overline{1}, -\tfrac{1}{4}) \\
&\sim c_{\overline{1},-\frac{1}{4}}(\varphi_1) \cdot H(\overline{0}, 0) + c_{\overline{1},-\frac{1}{4}}(\varphi_2) \cdot \left(\frac{105457575250}{169227} H(\overline{0}, 0) + H(\overline{0}, -1)\right) \\
&= \frac{1882717700}{169227} H(\overline{0}, 0) - \frac{1}{56} \cdot \left(\frac{105457575250}{169227} H(\overline{0}, 0) + H(\overline{0}, -1)\right) \\
&= -\frac{75}{28} H(\overline{0}, 0) - \frac{1}{56} H(\overline{0}, -1) \\
&= \frac{75}{28} \lambda - \frac{1}{56} D_{0,0} \ .
\end{aligned} \tag{19}$$

### 3.4.1 Writing the Hodge class in terms of Noether–Lefschetz divisors

We get interesting relations in $\mathrm{Pic}_{\mathbb{Q}}(\mathcal{F}_{2d})$ by writing the Hodge class $\lambda$ in terms of Noether–Lefschetz divisors. There is more than one way to do this, as there are infinitely many distinct Noether–Lefschetz divisors. We get a particularly interesting relation if we take as basis a simple set of Noether–Lefschetz divisors with low values of $\Delta$: see table 1.

Table 1: The Hodge relation for low values of $d$.

| $d$ | Hodge relation |
|---|---|
| 1 | $150\lambda \sim H(\overline{0}, -1) + 56 H(\overline{1}, -1/4)$ |
| 2 | $108\lambda \sim H(\overline{0}, -1) + 128 H(\overline{1}, -1/8) + 14 H(\overline{2}, -1/2)$ |
| 3 | $98\lambda \sim H(\overline{0}, -1) + 108 H(\overline{1}, -1/12) + 54 H(\overline{2}, -1/3) + 2 H(\overline{3}, -3/4)$ |
| 4 | $80\lambda \sim H(\overline{0}, -1) + 112 H(\overline{1}, -1/16) + 56 H(\overline{2}, -1/4) + 16 H(\overline{3}, -9/16)$ |

More examples can be found in [28]. Many of the coefficients in these relations have a geometric interpretation, counting curves with special properties; see [28, section 4.4.2] for an example.

## 3.5 The effective cone of $\mathcal{F}_{2d}$

Now that we have a good description of the Picard group $\mathrm{Pic}_{\mathbb{Q}}(\mathcal{F}_{2d})$, we may next try to understand the effective cone inside the Picard group.

There is a natural subcone of the effective cone, generated by the irreducible Noether–Lefschetz divisors $P_{\Delta,\delta}$:

**Definition 3.17.** The Noether–Lefschetz cone $\mathrm{Eff}^{\mathrm{NL}}(\mathcal{F}_{2d}) \subseteq \mathrm{Eff}(\mathcal{F}_{2d})$ is the cone generated by the set of all irreducible Noether–Lefschetz divisors $P_{\Delta,\delta}$.

Note that the reducible Noether–Lefschetz divisors $H(\gamma, n)$ are positive linear combinations of the irreducible ones, so they will all lie in this subcone.

First of all, it would be nice to understand the structure of the Noether–Lefschetz cone. In particular:

**Question 3.18.** Is the Noether–Lefschetz cone finitely generated? If so, can we give a list of generators and/or a list of bounding hyperplanes?

Apart from this, it would be nice to know how much we lose by restricting to this special subcone:

**Question 3.19.** Is the Noether–Lefschetz cone $\mathrm{Eff}^{\mathrm{NL}}(\mathcal{F}_{2d})$ equal to the effective cone $\mathrm{Eff}(\mathcal{F}_{2d})$ ?

Note that the equality $\mathrm{Pic}^{\mathrm{NL}}_{\mathbb{Q}}(\mathcal{F}_{2d}) = \mathrm{Pic}_{\mathbb{Q}}(\mathcal{F}_{2d})$, recently proved by [3], does not imply a positive answer to this last question: there could be effective divisors on $\mathcal{F}_{2d}$ that are linearly equivalent to some combination of irreducible Noether–Lefschetz divisors, but with some of the coefficients necessarily negative.

Also note that if $\mathrm{Eff}^{\mathrm{NL}}(\mathcal{F}_{2d}) \neq \mathrm{Eff}(\mathcal{F}_{2d})$, then there are modular forms on $\mathcal{F}_{2d}$ with a vanishing locus containing prime divisors that are not Noether–Lefschetz divisors.

We have done computer calculations to get an idea of the structure of the Noether–Lefschetz cone for many values of $d$; see [28, section 4.5] for details. The results suggest that it is generated by a relatively small number of irreducible Noether–Lefschetz divisors, all of small discriminant $\Delta$.

## 3.6 Deciding effectivity of the canonical class

In the cases for which we were able to compute the Noether–Lefschetz cone completely (up to $d = 32$), we may use those results to compute whether the canonical class $K^\circ$ is inside or outside the cone.

However, for the cases in the interesting region – say $d$ around 40, where we expect the Kodaira dimension of $\mathcal{F}_{2d}$ to change – it seemed not feasible to compute the Noether–Lefschetz cone in full detail. Fortunately, we do not need the full structure of the cone per se: we just need to know the relative position of the canonical class with respect to the cone, and if it is inside, we need to express the canonical class explicitly as a positive combination of (irreducible) Noether–Lefschetz divisors (as input for the calculation of the boundary coefficients).

It turns out that we may formulate this as a so-called linear programming problem. We want to write

$$K^\circ = \sum_{\Delta, \delta} t_{\Delta, \delta} [P_{\Delta, \delta}] \,, \tag{20}$$

where $K^\circ = 19\lambda - 1/2 \cdot [B]$ is the canonical class restricted to the open part of the moduli space, the $P_{\Delta, \delta}$ are irreducible Noether–Lefschetz divisors (see section 3.2), and the $t_{\Delta, \delta}$ are non-negative rational numbers.

**Remark 3.20.** We must restrict to a finite subset of the (infinite) set of irreducible Noether–Lefschetz divisors in order to get a finite problem. We may use the information gathered in our calculations of section 3.5 to guess which ones suffice to generate the full Noether–Lefschetz cone. This is of course not rigorous, and as a result we cannot conclude with certainty that a point is outside the Noether–Lefschetz cone. However, we do not even know that this cone equals the full effective cone of $\mathcal{F}_{2d}$ (see question 3.19), so this only adds to the uncertainty of an argument that was already incomplete. Moreover, the results of this procedure and their perfect agreement with the results of [12] suggest that in practice no information is lost at all.

So, pick a finite set of irreducible Noether–Lefschetz divisors $P_{\Delta,\delta}$, and introduce corresponding variables $t_{\Delta,\delta}$. Equation (20) is then a set of (linear) constraints for these variables (one constraint for every coordinate on the $\mathbb{Q}$-vector space $\mathrm{Pic}_{\mathbb{Q}}(\mathcal{F}_{2d})$). Existing programs for linear programming can solve these constraints, for non-negative values of the variables $t_{\Delta,\delta}$ – or assert that this is impossible.

We now discuss the implications of the results for the Kodaira dimension of $\mathcal{F}_{2d}$. We start with the lowest values of $d$, where $K^\circ$ is most negative.

**Theorem 3.21.** *If $1 \le d \le 15$, the moduli space $\mathcal{F}_{2d}$ has Kodaira dimension $-\infty$.*

*Proof.* Let us assume that the Kodaira dimension is not $-\infty$. Then there is some positive integer $m$ such that the multiple $mK$ is effective on $\overline{\mathcal{F}_{2d}}$. Restricting to the open part of the moduli space, we conclude that $mK^\circ$ is effective on $\mathcal{F}_{2d}$. It follows that the intersection product $mK^\circ \cdot \lambda^{18}$ is positive, because $\lambda$ is ample. Now, the number $K^\circ \cdot \lambda^{18} = \deg(K^\circ)$ is, up to a positive constant, given by the (vector-valued) Eisenstein series $-E_{\bar{0}}$ (seen as a function on $\mathrm{Pic}_{\mathbb{Q}}(\mathcal{F}_{2d})$) applied to $K^\circ$. However, our computations show that if $1 \le d \le 15$, then $K^\circ$ is not in the positive half-space determined by $E_{\bar{0}}$. This contradiction proves that $\kappa(\mathcal{F}_{2d}) = -\infty$. □

**Remark 3.22.** The use of intersection numbers on the quasi-projective variety $\mathcal{F}_{2d}$ may appear to be dubious. However, because the boundary components in the Satake compactification have very low dimension (only 0 and 1), and because $\lambda$ is ample even on this compactification, we may represent the class $\lambda^{18}$ – or in fact even $\lambda^2$ – on the compactification by a subvariety that is supported away from the boundary.

Most of these cases have been known for a long time by the work of Mukai ([20] [21] [22] [23] [24]), which uses the more explicit structure of the moduli space that is known in these cases. Our proof is simpler, using only coefficients of Eisenstein series and computation of the branch divisor. Moreover, the cases $d \in \{13, 14\}$ are new.

We now turn to the intermediate values of $d$ where $K^\circ$ is inside the positive half-space, but outside the Noether–Lefschetz cone.

**Theorem 3.23.** *Let $16 \le d \le 39$ or $d \in \{41, 44, 45, 47\}$. If the effective cone of $\mathcal{F}_{2d}$ is generated by irreducible Noether–Lefschetz divisors and our list of generators is complete (see questions 3.18 and 3.19), then $\kappa(\mathcal{F}_{2d}) = -\infty$.*

*Proof.* Solving the linear programming problem shows that for these values of $d$ the open part $K^\circ$ of the canonical class cannot be written as a non-negative combination of the supposedly generating Noether–Lefschetz divisors. By assumption these divisors indeed generate the effective cone, so we conclude that $K^\circ$ is not effective. Then $K = K^\circ - \Delta$ is definitely not effective, so $\kappa(\mathcal{F}_{2d}) = -\infty$. □

An unconditional proof that $\kappa(\mathcal{F}_{2d}) = -\infty$ in these cases thus needs a positive answer to question 3.19. It is suggestive that these values of $d$ (together with the cases $1 \le d \le 15$) are exactly the ones for which the alternative approach of [12], which aims to prove that $\kappa(\mathcal{F}_{2d}) \ge 0$, fails.

We may formulate this result positively as follows:

**Theorem 3.24.** *Let $16 \leq d \leq 39$ or $d \in \{41, 44, 45, 47\}$. Either $\kappa(\mathcal{F}_{2d}) = -\infty$ or there exists an irreducible codimension 1 subvariety of $\mathcal{F}_{2d}$ that is not a Noether–Lefschetz divisor.*

Finally, we turn to the cases where $K^\circ$ is inside the Noether–Lefschetz cone.

We can make a further distinction, by looking at the expression $K^\circ - \varepsilon\lambda$, where $\varepsilon$ is some nonnegative rational number. Recall that the canonical divisor $K$ is big if and only if $K - \varepsilon\lambda$ is effective for some positive number $\varepsilon$. So, if that is the case, then its restriction to the open part $\mathcal{F}_{2d}$, which is $K^\circ - \varepsilon\lambda$, must also be effective. Under the assumption that the answers to questions 3.19 and 3.18 are positive, we may compute whether this is possible. We simply extend the linear programming problem described by equation (20) by adding another variable $\varepsilon$, and changing the equations to

$$K^\circ - \varepsilon\lambda = \sum_{\Delta,\delta} t_{\Delta,\delta}[P_{\Delta,\delta}] \ . \tag{21}$$

We instruct the linear solver to minimise the solution with respect to the value of the variable $\varepsilon$.

The full results of this procedure, in the form of equations for $K^\circ - \varepsilon\lambda$ as a sum of irreducible Noether–Lefschetz divisors $P_{\Delta,\delta}$, can be found in [28]. We list in table 2 the most important part of the results: the minimal value of $\varepsilon$ among solutions of equation (21).

Table 2: The minimal value of $\varepsilon$ among solutions of equation (21).

| $d$ | 40 | 42 | 43 | 46 | 48 | 49 | 50 | 51 | 52 | 53 | 54 | 55 | 56 |
|---|---|---|---|---|---|---|---|---|---|---|---|---|---|
| $\varepsilon$ | 0 | 0 | 0 | 1 | 0 | 0 | 1 |  | 1 |  | 1 | 0 | 0 |

Now, let us consider what these data mean for the Kodaira dimension of $\mathcal{F}_{2d}$.

**Theorem 3.25.** *Let $d \in \{40, 42, 43, 48, 49, 55, 56\}$. If the effective cone of $\mathcal{F}_{2d}$ is generated by irreducible Noether–Lefschetz divisors and our list of generators is complete, then the Kodaira dimension of $\mathcal{F}_{2d}$ satisfies $\kappa(\mathcal{F}_{2d}) < 19$.*

*Proof.* We see from table 2 that for these values of $d$, equation (21) can only be solved for $\varepsilon = 0$; therefore, by the assumption on the effective cone, $K^\circ - \varepsilon\lambda$ cannot be effective for positive $\varepsilon$. In such a case $K$ cannot be big, so the Kodaira dimension of $\mathcal{F}_{2d}$ is less than 19. □

Again, we may formulate this result positively as follows:

**Theorem 3.26.** *Let $d \in \{40, 42, 43, 48, 49, 55, 56\}$. Either $\kappa(\mathcal{F}_{2d}) < 19$, or there exists an irreducible codimension 1 subvariety of $\mathcal{F}_{2d}$ that is not a Noether–Lefschetz divisor.*

In fact, for these values of $d$ we can also prove unconditionally that $\kappa(\mathcal{F}_{2d}) \geq 0$, but for that we need to consider the boundary of the moduli space; see theorem 4.14.

The agreement with the results of [12] is again striking: these $d$ are exactly the values for which their method only proves that $\kappa(\mathcal{F}_{2d}) \geq 0$.

For the other values of $d$ (i.e., the ones giving a solution to (21) with positive $\varepsilon$) the divisor $K - \varepsilon\lambda$ on $\overline{\mathcal{F}_{2d}}$ stands a chance of being effective, which would mean that $\kappa(\mathcal{F}_{2d}) = 19$. However, in order to prove this we need to take the boundary of the moduli space into account. We would like to take the expression of $K° - \varepsilon\lambda$ as a combination of Heegner divisors and extend it to a valid relation in the Picard group of $\overline{\mathcal{F}_{2d}}$. In the next section, we show how to do this, by computing boundary coefficients of relations among Noether–Lefschetz divisors (see theorem 4.8); in theorem 4.12 we apply this to prove that $\mathcal{F}_{2d}$ is of general type for the relevant values of $d$.

# 4 The moduli space of polarised K3 surfaces: boundary

## 4.1 Satake compactification

From now on we will write $L = L_{2d} = \langle -2d \rangle \oplus 2U \oplus 2E_8(-1)$ for brevity.

The Satake compactification, also called Baily–Borel compactification,

$$\mathcal{F}_{2d}^* = \mathrm{Proj}\Big(\bigoplus_k M(k, \tilde{\mathrm{O}}^+(L))\Big) \tag{22}$$

adds to $\mathcal{F}_{2d}$ a finite number of 0- and 1-dimensional components called cusps. We now briefly review the properties of the Satake compactification of $\mathcal{F}_{2d}$; for more details, see section 5.1 of [28].

The 0-cusps correspond to isotropic elements of the discriminant group $D_L$ (i.e., $\gamma \in D_L$ such that $\gamma^2/2 = 0 \in \mathbb{Q}/\mathbb{Z}$). There are few of these (at most four 0-cusps for $d \leq 61$) and they are easy to compute.

The 1-cusps correspond to isotropic planes $I \subset L \otimes \mathbb{Q}$ up to the action of $\tilde{\mathrm{O}}^+(L)$. The set of 1-cusps is large (it grows as $d^8$) and hard to compute exactly.

However, we may associate a definite lattice $K(F)$ of rank 17 to every 1-cusp $F$ (to wit, if $F$ is represented by the isotropic plane $I$, $K(F)$ is the subquotient $I^\perp/I$), and for our purposes it suffices to know which definite lattices arise in this way.

These definite lattices $K(F)$ have a discrete invariant $N$ called the imprimitivity. Roughly speaking, $N$ measures how much of the discriminant group $\mathbb{Z}/2d\mathbb{Z}$ of $L$ is no longer present in the discriminant group of $K(F)$.

**Definition 4.1** ([32])**.** Let $I$ be an isotropic plane in $L$. Define $H_I = (I_{L^\vee}^\perp)_{L^\vee}^\perp / I$; if $H_I \cong \mathbb{Z}/N\mathbb{Z}$, then the positive integer $N$ is by definition the imprimitivity of $I$. If the 1-cusp $F$ is represented by $I$, then we also write $H_F$ for $H_I$ and call $N$ the imprimitivity of $F$.

The lattice genus of the definite lattice $K(F)$ is not the same for all 1-cusps $F$. However, it only depends on the imprimitivity invariant $N$ of the cusp:

**Proposition 4.2** ([32, Lemma 5.1.3])**.** *If the cusp $F$ has imprimitivity $N$, then the discriminant module of $K(F)$ is isomorphic to $\mathbb{Z}/2m\mathbb{Z}$, where $m = d/N^2$.*

We will need the following explicit version of the above statement, giving the relation between the groups $\mathbb{Z}/2d\mathbb{Z}$ and $\mathbb{Z}/2m\mathbb{Z}$. Suppose that the 1-cusp $F$ is represented by the isotropic plane $I$. The discriminant group $D_K \cong \mathbb{Z}/2m\mathbb{Z}$ of the subquotient lattice $K(F) = I^\perp/I$ of $L$ is naturally a subquotient of the discriminant group $D_L$: it is $H_F^\perp/H_F$ (see below).

**Definition 4.3.** The subgroup $H_F^\perp \subseteq D_L$ is the subgroup orthogonal to $H_F$:

$$H_F^\perp = \{\gamma \in D_L : \text{for all } \delta \in H_F : (\gamma, \delta) = 0\} \ . \tag{23}$$

We define $p : H_F^\perp \to D_K$ to be the surjective map derived from the isomorphism $D_K \cong H_F^\perp/H_F$.

For example, if the cusp $F$ has imprimitivity $N = 1$ (this is the case for the standard cusp, for instance), then $H_F = \{\bar{0}\}$, so $H_F^\perp = D_L$, and the map $p : D_L \to D_K$ is an isomorphism of discriminant groups.

Now, proposition 4.2 implies that if $F$ has imprimitivity $N$, then $K(F) \in \mathcal{G}(K(I_N)) = \mathcal{G}(\langle -2m \rangle \oplus 2E_8(-1))$.

## 4.2 Toroidal compactification

We will briefly describe toroidal compactifications $\overline{\mathcal{F}_{2d}}$ of our locally symmetric domain $\mathcal{F}_{2d}$. We follow the notation and description of [15, section 5.3]; see also [28, section 5.2], [12] and the book [1].

Given a cusp $F$, we have the stabiliser (parabolic subgroup) $N(F) \subset \tilde{O}^+(L)_\mathbb{R}$, the unipotent radical $W(F) \subset N(F)$, and the centre $U(F) \subseteq W(F)$. The partial compactification of $\mathcal{D}_{2d}$ at $F$ is taken inside the larger space $\mathcal{D}_L(F)$. This space $\mathcal{D}_L(F)$ can be abstractly defined as

$$\mathcal{D}_L(F) = U(F)_\mathbb{C} \mathcal{D}_{2d} \ , \tag{24}$$

where $U(F)_\mathbb{C}$ acts on the period domain $\mathcal{D}_{2d}$ within the larger space (the so-called compact dual)

$$\check{\mathcal{D}}_{2d} = \{\mathbb{C}z : (z, z) = 0\} \subset \mathbb{P}(L \otimes \mathbb{C}) \ . \tag{25}$$

The space $\mathcal{D}_L(F)$ has a product decomposition

$$\mathcal{D}_L(F) \cong F \times V(F) \times U(F)_\mathbb{C} \ , \tag{26}$$

where $V(F) = W(F)/U(F)$.

### 4.2.1 1-cusps

Suppose that the cusp $F$ is a 1-cusp corresponding to an isotropic plane in $L$. The stabiliser $N(F)$, unipotent radical $W(F)$ and its centre $U(F)$ can be described explicitly with respect to a particular basis of $L$: see [12, section 2.3]. The subgroup $U(F)$ is 1-dimensional in this case.

The decomposition (26) in this case looks like

$$\mathcal{D}_L(F) \cong \mathbb{C} \times \mathbb{C}^{17} \times \mathbb{H} \ ; \tag{27}$$

let us write $s$ for the coordinate on the first factor $\mathbb{C}$.

The group $U(F)$ is 1-dimensional, so the torus $U(F)_\mathbb{C}/U(F)_\mathbb{Z}$ is just $\mathbb{C}^\times$. A calculation shows that the element of $U(F)$ parametrised by $x \in \mathbb{R}$ acts on $\mathcal{D}_L(F)$ by increasing $s$ by $Nx$ (and fixing the other coordinates) – see [12, proposition 2.26]. Therefore, we choose as a coordinate on the torus $U(F)_{\mathbb{C}}/U(F)_{\mathbb{Z}} \cong \mathbb{C}^{\times}$ the function $u = \exp(2\pi \mathrm{i} s/N)$; the compactification then adds the point $u = 0$.

Because the group $U(F)$ is 1-dimensional, the real cone $C(F)$ is just $\mathbb{R}_+$, and for such a trivial cone there is only a single choice of fan. Therefore, as far as the toroidal boundary over the 1-cusps is concerned, we do not need to make any choices.

### 4.2.2 0-cusps

Suppose that the cusp $F$ is a 0-cusp, corresponding to a primitive isotropic vector $z \in L$. Writing $M = z^{\perp}/\mathbb{Z}z$ as before, the stabiliser is given by the semi-direct product $N(F) = \tilde{\mathrm{O}}^+(M) \rtimes E_z(L)$, where $E_z(L)$ is the group of Eichler transvections associated to $z$:

**Definition 4.4.** Let $L$, $z$ and $M$ be as above. For any $m \in M$, the Eichler transvection associated to $m$ is the map $L \to L$ given by

$$x \mapsto x - (x,z)m + (x,m)z - m^2/2 \cdot (x,z)z \ . \tag{28}$$

The group of all such transformations is denoted by $E_z(L)$.

In fact, this group of Eichler transvections is isomorphic to the additive group of the lattice $M$.

Furthermore, the unipotent radical $W(F)$ is the subgroup $E_z(L) \cong M$; as this is an abelian group, its centre $U(F)$ coincides with $W(F)$.

Now, looking at equation (26), we see that in this case $F \cong \{\mathrm{pt}\}$ (as $F$ is a 0-cusp), and $V(F) = W(F)/U(F) = \{0\}$, so the space $\mathcal{D}_L(F)$ is actually isomorphic to $U(F)_{\mathbb{C}} \cong M \otimes \mathbb{C}$. In fact, the inclusion $\mathcal{D}_{2d} \subseteq \mathcal{D}_L(F)$ is just the tube domain realisation associated to $z$ (see section 5.2.1 of [28]), so under the isomorphism $\mathcal{D}_L(F) \cong M \otimes \mathbb{C}$, the period domain $\mathcal{D}_{2d}$ is given by the set $M \otimes \mathbb{R} \oplus iC(F)$, where $C(F)$ is a real cone in $M \otimes \mathbb{R}$. In fact, $C(F)$ is just the positive cone in $M \otimes \mathbb{R}$, given by $C(F) = \{v \in M \otimes \mathbb{R} : v^2 > 0, (z',v) > 0\}$. (Recall that $z'$ is an element of $L \otimes \mathbb{Q}$ such that $(z,z') = 1$; here, it just serves to pick out one of the two connected components of the set of vectors of positive norm.)

The torus $T(F)$ is $U(F)_{\mathbb{C}}/U(F)_{\mathbb{Z}} \cong (M \otimes \mathbb{C})/M \cong (\mathbb{C}^{\times})^{19}$. The toroidal compactifications complete this torus to a toric variety, but this now essentially depends on a choice of fan $\Sigma(F)$ (i.e., cone decomposition) of the rational closure of the cone $C(F)$, a real cone in the 19-dimensional vector space $M \otimes \mathbb{R}$.

The components of the boundary divisor over a cusp in general correspond to the rays in the cone decomposition of the rational closure of $C(F)$, up to action by the orthogonal group $\tilde{\mathrm{O}}^+(M)$ and identification in the final gluing procedure. In the case of a 0-cusp, this set of rays depends on the choice of fan $\Sigma(F)$: at the very least, it will include the boundary rays of the rational closure of $C(F)$, but there may be more (internal) rays.

Note that the boundary rays of the rational closure of the positive cone $C(F) \subset M \otimes \mathbb{R}$ are exactly the rays through the isotropic vectors of $M$. Such a vector, taken together with the isotropic vector $z \in L$ representing the 0-cusp under consideration, gives an isotropic plane in $L$, which in turns represents a 1-cusp of $\mathcal{F}_{2d}$. Moreover, the component of the boundary divisor over the 0-cusp that corresponds to this ray is identified by the gluing procedure with the boundary component over this 1-cusp (as described in section 4.2.2).

There is one somewhat natural choice of fan in this case: the perfect cone decomposition (see [1] for details). It is minimal, in the sense that the set of rays in this decomposition consists of only the boundary rays. This is useful: by the above paragraph, we see that if we choose this fan for every 0-cusp, the only boundary components that we get are the ones over 1-cusps. The main disadvantage of the perfect cone decomposition is that it gives a compactification that may have bad (i.e., non-canonical) singularities. A further subdivision of the decomposition is necessary to get rid of these singularities ([12, section 2] proves that that is possible), but this reintroduces additional boundary components that are harder to control.

One feature of our case simplifies the situation significantly: because the cone decompositions over the 1-cusps are unique, the compatibility conditions between the fans associated to a 1-cusp and 0-cusps in its closure are trivially satisfied.

### 4.2.3 Sufficient conditions for cuspidality

For our applications we want modular forms on $\mathcal{F}_{2d}$ that are cusp forms in a very strong sense: for any toroidal compactification, we want the modular form to vanish on every component of the boundary divisor. In fact, because the vanishing order need not be an integer number, we should strengthen this by demanding that the vanishing order is at least 1.

**Lemma 4.5** ([13]). *If the modular form $\Psi$ on $\mathcal{F}_{2d}$ vanishes at all cusps, then it vanishes to order at least 1 on every component of the boundary divisor in any toroidal compactification.*

This fact is used implicitly in [12], without proof. The essential point is that the vanishing order of a modular form on $\mathcal{F}_{2d}$ at a component of the boundary divisor is an integer number; this is proved in [14, proposition 2.1] (also see [13]). The proof relies on the fact that the arithmetic group $\tilde{O}^+(L_{2d})$ has only one non-trivial character (the determinant); see [14, corollary 1.8].

Moreover, it is enough to have cuspidality at the 1-cusps:

**Proposition 4.6.** *If the modular form $\Psi$ on $\mathcal{F}_{2d}$ vanishes at all 1-cusps, then it vanishes to order at least 1 on every component of the boundary divisor in any toroidal compactification.*

*Proof.* Every 0-cusp occurs as a limit point of at least one 1-cusp, so $\Psi$ vanishes at every 0-cusp as well by continuity. Lemma 4.5 gives the desired result. □

## 4.3 Extending divisor relations to the compactification

To complete a given relation in $\mathrm{Pic}_{\mathbb{Q}}(\mathcal{F}_{2d})$ to a relation in $\mathrm{Pic}_{\mathbb{Q}}(\overline{\mathcal{F}_{2d}})$, we will proceed as follows. Recall that the map $c_{\gamma,n} \mapsto H(\gamma, n)$ sends identities among coefficients of modular forms to linear equivalences of divisors. The function – in fact, a meromorphic modular form – on $\mathcal{F}_{2d}$ exhibiting this linear equivalence can be described quite explicitly, by the work of Borcherds. We will determine the vanishing order of this function at the cusps, and consequently compute the boundary terms of its divisor. This will give a relation in $\mathrm{Pic}_{\mathbb{Q}}(\overline{\mathcal{F}_{2d}})$.

**Remark 4.7.** In view of the discussion of section 4.2.2, for this to make sense we need to pick a specific toroidal compactification, because that choice determines the structure of the boundary divisors over the 0-cusps.

We pick the toroidal compactification determined by the perfect cone decomposition. We have seen in section 4.2.2 that in that case the irreducible components of the boundary divisor are easy to describe: there is one component for every 1-cusp, and there are no others. Let us write $\Delta_F$ for the component corresponding to the 1-cusp $F$.

Now, we can formulate the boundary behaviour of the modular form $\Psi$.

**Theorem 4.8.** *Let $\sum_{\gamma,n} a_{\gamma,n} H(\gamma, n) \sim 0$ be a linear equivalence of Noether–Lefschetz divisors on $\mathcal{F}_{2d}$. Then the following linear equivalence holds on $\overline{\mathcal{F}_{2d}}$ (the toroidal compactification of $\mathcal{F}_{2d}$ with the perfect cone decomposition):*

$$\sum_{\gamma,n} a_{\gamma,n} H(\gamma, n) + \sum_{F \in S_1} \sum_{\gamma,n} a_{\gamma,n} c(\gamma, n, N_F, K(F)) \Delta_F \sim 0 \ . \quad (29)$$

*Here $F$ ranges over the 1-cusps $S_1$ of $\mathcal{F}_{2d}$, $N_F$ is the imprimitivity of the cusp $F$ (see section 4.1), and $K(F)$ is the negative definite lattice of rank 17 associated to the cusp $F$ (also explained in section 4.1). The function $c(\gamma, n, N_F, K(F))$ calculating the contribution of a given Heegner divisor $H(\gamma, n)$ at the cusp of imprimitivity $N_F$ having definite lattice $K(F)$ is given by*

$$c(\gamma, n, N, K) = \begin{cases} N/24 \cdot (E_2 \boldsymbol{\Theta}_K)(p(\gamma), n) & \text{if } \gamma \in H_F^\perp \\ 0 & \text{otherwise} \end{cases} , \quad (30)$$

*where $E_2$ is the usual Eisenstein series, $\boldsymbol{\Theta}_K$ is the vector-valued theta series of the lattice $K$, and the subgroup $H_F^\perp \subseteq D_L$ and the map $p : H_F^\perp \to D_K$ are defined in 4.3.*

The proof of this theorem can be found in section 5.3.4 of [28]; it amounts to an analysis of the limit behaviour of Borcherds' product formula (12) as we approach the boundary of $\mathcal{F}_{2d}$.

We look at an example application of this theorem in the next section.

## 4.4 Example: completing the Hodge relation for $d = 1$

As an example, we take the Hodge relation in $\text{Pic}_\mathbb{Q}(\mathcal{F}_{2\cdot 1})$:

$$150 \lambda \sim H(\bar{0}, -1) + 56 H(\bar{1}, -1/4) \ . \quad (31)$$

For $d = 1$, there is one 0-cusp, and there are four 1-cusps, $F_\beta, F_\gamma, F_\zeta, F_\eta$ (see the example 5.1.16 of [28]). Let us write $\Delta = \Delta_\beta + \Delta_\gamma + \Delta_\zeta + \Delta_\eta$ for the corresponding decomposition of the boundary divisor.

These four 1-cusps give only two distinct vector-valued theta series: $\boldsymbol{\Theta}_{K_\beta} = \boldsymbol{\Theta}_{K_\gamma}$ and $\boldsymbol{\Theta}_{K_\zeta} = \boldsymbol{\Theta}_{K_\eta}$. Applying theorem 4.8 to this relation, we get

$$\begin{aligned}
0 \sim\ & 150 H(\bar{0}, 0) + H(\bar{0}, -1) + 56 H(\bar{1}, -1/4) \\
& + (150 c(\bar{0}, 0, 1, K_\beta) + c(\bar{0}, 1, 1, K_\beta) + 56 c(\bar{1}, -1/4, 1, K_\beta))(\Delta_\beta + \Delta_\gamma) \\
& + (150 c(\bar{0}, 0, 1, K_\zeta) + c(\bar{0}, 1, 1, K_\zeta) + 56 c(\bar{1}, -1/4, 1, K_\zeta))(\Delta_\zeta + \Delta_\eta) \\
\sim\ & 150 H(\bar{0}, 0) + H(\bar{0}, -1) + 56 H(\bar{1}, -1/4) + 30 (\Delta_\beta + \Delta_\gamma) + 18 (\Delta_\zeta + \Delta_\eta) \ .
\end{aligned} \quad (32)$$

We thus see that

$$150\,\lambda \sim H(\overline{0},-1) + 56\,H(\overline{1},-1/4) + 30\,(\Delta_\beta + \Delta_\gamma) + 18\,(\Delta_\zeta + \Delta_\eta)\,. \tag{33}$$

Computations such as this one are possible in principle for every $d$. However, the number of 1-cusps increases very rapidly with $d$ (as $d^8$), and hence so does the size of the resulting equation.

## 4.5 Theta ghosts

We would like to use theorem 4.8 to complete relations in $\mathrm{Pic}(\mathcal{F}_{2d})$ to the boundary, in particular to prove that given modular forms on $\mathcal{F}_{2d}$ are cusp forms.

Looking at equations (29) and (30), this means that we have to enumerate the 1-cusps of $\mathcal{F}_{2d}$ and for each of them compute the associated definite lattice $K$ and its vector-valued theta series $\boldsymbol{\Theta}_K$. We know from proposition 4.2 that the lattice genus of $K$ is known (although dependent on the imprimitivity $N$ of $F$): it is $\mathcal{G}(2E_8(-1) \oplus \langle -2m \rangle)$ where $m = d/N^2$.

Let us rephrase our task in a slightly more general setting. Let $K_0$ be a non-degenerate definite lattice. We would like to determine what vector-valued modular forms occur as theta series of a definite lattice of genus $\mathcal{G}(K_0)$. (For an introduction to vector-valued modular forms, see section 2.1; for vector-valued theta series (these count the number of vectors in the dual lattice of given length and discriminant class), see 2.1.1.) In other words, we want to identify the image of the map

$$\boldsymbol{\Theta}\,:\,\mathcal{G}(K_0) \to M_{\mathrm{rank}(K_0)/2}(K_0)\,:\,K \mapsto \boldsymbol{\Theta}_K\,. \tag{34}$$

In our application we have $K_0 = 2E_8 \oplus \langle 2m \rangle$, where $m = d/N^2$ is directly related to the polarisation degree of the K3 surface (in particular there is always the case $N = 1, m = d$). As $m$ increases, computing the image of $\boldsymbol{\Theta}$ is computationally too hard to approach directly, i.e., by enumerating the lattices in the genus $\mathcal{G}(K_0)$ and computing the theta series of each lattice. This is already clear from the size of the set $\mathcal{G}(K_0)$, which grows rapidly (as $m^8$).

We will take another approach: we use some properties shared by all theta series to define a finite superset of the image of $\boldsymbol{\Theta}$. For some reasonably low values of $d$, we can compute this superset explicitly. More importantly, for all $d$ in the range that interests us we can compute a lower bound for the boundary coefficients of any given relation in $\mathrm{Pic}(\mathcal{F}_{2d})$ by solving a linear programming problem that minimises a linear expression for these boundary coefficients over the superset.

**Definition 4.9.** Let $K$ be a definite lattice (let us say positive definite). A vector-valued modular form $\Psi = \sum_{\gamma,n} \Psi(\gamma,n)\,q^n \mathbf{e}_\gamma \in M_{\mathrm{rank}\,K/2}(K)$ is a theta ghost if

(i) it is an almost-cusp form (see definition 2.7);

(ii) $\Psi(\gamma,n) \in \mathbb{N}$ for all $\gamma \in D_K$ and $n \geq 0$;

(iii) $\Psi(\gamma,n) = \Psi(-\gamma,n)$ for all $\gamma \in D_K$ and $n \geq 0$;

(iv) $\Psi(\gamma,n) \in 2\mathbb{N}$ for all $\gamma \in D_K$ such that $\gamma = -\gamma$ and for all $n > 0$;

(v) $\Psi(0,0) = 1$.

We see that if $K$ is a definite lattice, then its vector-valued theta series $\boldsymbol{\Theta}_K$ is a theta ghost.

We would like to compute the set of theta ghosts. Using the method outlined in section 3.4, we may compute a Heegner basis $c^\vee_{\gamma_i,n_i}$ of the space $M_{17/2}(K_0)$ and write any other coefficient function $c_{\gamma,n}$ as an explicit (rational) linear combination of these basis vectors. Write $\Psi = \sum_i m_i c^\vee_{\gamma_i,n_i}$ for a general theta ghost. For every choice of $\gamma$ and $n$, we rewrite condition (ii) on the coefficient $\Psi(\gamma, n)$ as a condition on the numbers $m_i$. This gives a restriction of the form

$$\sum_i a_i m_i \in \mathbb{N} , \qquad (35)$$

where the $a_i$ are rational numbers. Clearing denominators, we see that this is a problem of linear integer programming, and our task is to enumerate all solutions.

**Remark 4.10.** If we are given a modular form in a numerical way, where the calculation of every coefficient requires some amount of work, then it is impossible to verify the infinite set of conditions (ii) directly.

However, examples show that a finite number of these conditions (ii) seem to suffice: if a modular form satisfies condition (ii) for this finite set of coefficients $(\gamma, n)$, then apparently it is fulfilled for all other coefficients as well (as far as we could check). We will call such modular forms *apparent* theta ghosts. Note that this situation is analogous to the question of finite generation of the Noether–Lefschetz cone: see question 3.18.

Although the sets of theta ghosts and apparent theta ghosts apparently coincide, we do not know for sure that they do, and we will not use that fact. We do know that the set of theta ghosts is contained in the set of apparent theta ghosts, and thus the set of theta series is contained in the set of apparent theta ghosts; that last inclusion is enough for our purposes.

By solving the integer linear programming problem, we have counted the number of apparent theta ghosts for $d$ at most 10: see table 3.

Table 3: The number $G$ of apparent theta ghosts for given $d$.

| $d$ | 1 | 2 | 3 | 4 | 5 | 6 | 7 | 8 | 9 | 10 |
|---|---|---|---|---|---|---|---|---|---|---|
| $G$ | 3 | 35 | 11 | 107 | 58 | 164 | 483 | 1344 | 196 | 4887 |

## 4.6 Application to cuspidality of modular forms on $\mathcal{F}_{2d}$

We show how to use theta ghosts to prove that a given modular form on $\mathcal{F}_{2d}$ is a cusp form.

Given any modular form on $\mathcal{F}_{2d}$ that has vanishing locus supported on Noether–Lefschetz divisors, we know by Bruinier's result that it arises from Borcherds' construction (see section 3.3). Our theorem 4.8 then describes the vanishing order of the modular form at every 1-cusp (or more precisely: at the components over all the 1-cusps of any toroidal compactification). Because the number of 1-cusps increases rapidly with the polarisation degree $2d$, we use our idea of theta ghosts to bound these vanishing orders.

We recall here the contents of theorem 4.8, reinterpreted in terms of modular forms instead of relations in $\mathrm{Pic}_{\mathbb{Q}}(\mathcal{F}_{2d})$. Let $\Psi$ be any modular form on $\mathcal{F}_{2d}$ with vanishing locus supported on Noether–Lefschetz divisors; we may assume without loss of generality that $\Psi$ is associated to the linear relation $\sum_{\gamma,n} a_{\gamma,n} c_{\gamma,n} = 0$ of coefficients functions $c_{\gamma,n}$ on $AC(d)$, the space of almost cusp forms associated to the lattice $L_{2d}$. Then the vanishing order of $\Psi$ at the boundary components associated to the 1-cusp $F$ is given by

$$\sum_{\substack{\gamma,n \\ \gamma \in H_F^{\perp}}} a_{\gamma,n} \cdot N/24 \cdot (E_2 \Theta_K)(p(\gamma), n) , \tag{36}$$

where $N$ is the imprimitivity of the cusp $F$ (see section 4.1), $K$ is the definite lattice of rank 17 associated to the cusp $F$ (see section 4.1), $\Theta_K$ is the vector-valued theta series of the lattice $K$ (see section 2.1.1), and $E_2$ is the usual Eisenstein series.

For a fixed modular form $\Psi$, this vanishing order is thus a linear expression in the coefficients $\Theta_K(\gamma, n)$.

For all 1-cusps $F$ of the same imprimitivity $N$, the vector-valued theta series $\Theta_K$ is a member of the space $AC_{17/2}(m)$ almost cusp forms of weight $17/2$ associated to the lattice $\mathbb{Z}/2m\mathbb{Z}$, where $m = d/N^2$. We may compute a basis of that space (see section 3.4) and express any coefficient $c_{\gamma,n}$ in terms of that basis. Let us write $\lambda_i$ for the unknown coefficients in the expression for $\Theta_K$ in terms of the basis. Now, the fact that the coefficients $\Theta_K(\gamma, n)$ are natural numbers gives us an infinite set of inequalities and integrality constraints for the numbers $\lambda_i$.

Summarising, we have a linear expression of the vanishing order of $\Psi$ – viewed as a function of the cusp $F$ – in terms of variables $\lambda_i$ and we have linear inequalities and integrality constraints for these $\lambda_i$. We may apply integer linear programming techniques to minimise the expression (36) with respect to the constraints. This gives us a lower bound for the vanishing order of $\Psi$ among all the 1-cusps with imprimitivity $N$.

Repeating this for all possible imprimitivities $N$ (i.e., the positive integers with $N^2$ dividing $d$), and taking the lowest outcome, we get a lower bound for the vanishing order of $\Psi$ among all the 1-cusps. If this is a positive number, then we conclude that $\Psi$ is a cusp form.

Applying this procedure to the modular forms we constructed on $\mathcal{F}_{2d}$ we get lower bounds for the vanishing orders. We find that all the lower bounds are positive, so we conclude that all the corresponding modular forms are cusp forms.

**Remark 4.11.** Because there might be theta ghosts that do not come from definite lattices, this method cannot prove that a form is not a cusp form.

### 4.7 Application to the Kodaira dimension of $\mathcal{F}_{2d}$

**Theorem 4.12.** *If $d \in \{46, 50, 52, 54\}$, then $\mathcal{F}_{2d}$ is of general type (i.e., the Kodaira dimension is $\kappa(\mathcal{F}_{2d}) = 19$).*

*Proof.* By the results of section 3.6, for these values of $d$ we can write $K^\circ - \varepsilon\lambda = (19-\varepsilon)\lambda - 1/2 \cdot [B]$ as a positive linear combination of (irreducible) Noether–Lefschetz divisors for some positive value of $\varepsilon$:

$$(19 - \varepsilon)\lambda - 1/2 \cdot [B] = \sum_{\Delta,\delta} t_{\Delta,\delta}[P_{\Delta,\delta}] \ . \tag{37}$$

Translating, this means that there is a modular form $\Psi$ on $\mathcal{F}_{2d}$ of weight $19 - \varepsilon < 19$ that vanishes on the ramification divisor (which is exactly the pullback of $1/2 \cdot B$). By [12, theorem 1.1], if we can prove that $\Psi$ is a cusp form, then we know that $\mathcal{F}_{2d}$ is of general type.

First of all, we rewrite the irreducible Noether–Lefschetz divisors $P_{\Delta,\delta}$ in terms of the reducible divisors $H(\gamma,n)$ using the triangular relations of [28, section 4.2.2].

Next, we apply theorem 4.8 to relation (37), completing it to a relation on some toroidal compactification. In theory, this gives us the boundary coefficient of (37) at every 1-cusp. The set of 1-cusps is very large, though, so we cannot compute all of these boundary coefficients explicitly.

We apply the method of theta ghosts (see section 4.5) to compute a lower bound for the boundary coefficient of $\Psi$ at all the 1-cusps: we take the expression produced by theorem 4.8 for the vanishing order in terms of theta coefficients, and use linear programming to get a lower bound for this expression. From the results of this procedure, we see that for the values of $d$ under consideration the lower bound is at least 1. This means that the modular form $\Psi$ vanishes at every 1-cusp.

By proposition 4.6, this is enough to conclude that $\Psi$ is indeed a cusp form, and we are done. □

**Remark 4.13.** Note that in the proof we use two distinct toroidal compactifications of $\mathcal{F}_{2d}$: theorem 4.8 (which completes relation on $\mathcal{F}_{2d}$ to $\overline{\mathcal{F}_{2d}}$) applies to the toroidal compactification associated to the perfect cone decomposition; in the proof of [12, theorem 1.1] on the other hand it is essential to use a toroidal compactification with canonical singularities. This is not a problem: the only thing we need from theorem 4.8 is the vanishing order of $\Psi$ at the 1-cusps, and this is independent of the choice of toroidal compactification.

**Theorem 4.14.** *If $d \in \{40, 42, 43, 48, 49, 55, 56\}$, then $\kappa(\mathcal{F}_{2d}) \geq 0$.*

*Proof.* The only difference with the theorem above is that we now have a modular form of weight exactly 19 (because $\varepsilon = 0$ in these cases). This means we can apply the second case of [12, theorem 1.1], and conclude that the Kodaira dimension of $\mathcal{F}_{2d}$ is non-negative. (Again, the vanishing order at the standard 1-cusp of the form we find equals 15 in all of these cases.) □

Combining this result with theorem 3.25, we get the following.

**Theorem 4.15.** *Let $d \in \{40, 42, 43, 48, 49, 55, 56\}$. If the effective cone of $\mathcal{F}_{2d}$ is generated by irreducible Noether–Lefschetz divisors and our list of generators is complete (see questions 3.18 and 3.19), then we have intermediate Kodaira dimension: $0 \leq \kappa(\mathcal{F}_{2d}) < 19$.*